\input amstex
\input epsf
\documentstyle{amsppt}\nologo\footline={}

\def\ad{\mathop{\text{\rm ad}}}
\def\arccosh{\mathop{\text{\rm arccosh}}}
\def\d{\mathop{\text{\rm d}}}
\def\dist{\mathop{\text{\rm dist}}}
\def\Gr{\mathop{\text{\rm Gr}}}
\def\Lin{\mathop{\text{\rm Lin}}}
\def\ker{\mathop{\text{\rm ker}}}
\def\GL{\mathop{\text{\rm GL}}}
\def\T{\mathop{\text{\rm T}}}
\def\U{\mathop{\text{\rm U}}}
\def\S{\mathop{\text{\rm S}}}
\def\Re{\mathop{\text{\rm Re}}}
\def\tr{\mathop{\text{\rm tr}}}
\def\ch{\mathop{\text{\rm char}}}
\def\Spin{\mathop{\text{\rm Spin}}}
\def\Sp{\mathop{\text{\rm Sp}}}
\def\ta{\mathop{\text{\rm ta}}}
\def\Area{\mathop{\text{\rm Area}}}

\hsize450pt

\topmatter\title Grassmannians and conformal structure on
absolutes\endtitle\author Sasha Anan$'$in, Eduardo C.~Bento
Gon\c calves, and Carlos H.~Grossi\endauthor\thanks First author
partially supported by the Institut des Hautes \'Etudes Scientifiques
(IHES).\endthanks\thanks Second author supported by FAPESP, project
No.~07/56989-0.\endthanks\thanks Third author partially supported by the
Max-Planck-Gesellschaft.
\endthanks\address Departamento de Matem\'atica, ICMC, Universidade de S\~ao Paulo,
Caixa Postal 668,\newline13560-970--S\~ao Carlos--SP, Brasil\endaddress\email
sasha\_a\@icmc.usp.br\endemail\address Independent Scholar,
Campinas--SP, Brasil\endaddress\email
EduardoCBG\@gmail.com\endemail\address Departamento de Matem\'atica, ICMC, Universidade de S\~ao
Paulo, Caixa Postal 668,\newline13560-970--S\~ao Carlos--SP, Brasil\endaddress\email
grossi\@icmc.usp.br\endemail\subjclass 53A20 (53A35,
51M10)\endsubjclass\abstract We study grassmannians associated with a linear space with a
nondegenerate hermitian form. The geometry of these grassmannians allows us to explain the relation
between a (pseudo-)riemannian projective geometry and the conformal structure on its ideal
boundary (absolute). Such relation encompasses, for instance, the usual conformal structure on the
absolute of real hyperbolic space, the usual conformal structure on the absolute of de Sitter space,
the conformal contact structure on the absolute of complex hyperbolic space, and the causal structure
on the absolute of anti-de Sitter space.
\endabstract\endtopmatter\document

\vskip-10pt

\hfill
{\it Dedicated to the memory of Waldyr Rodrigues Jr.}

\vskip30pt

\centerline{\bf1. Introduction}

\medskip

We study the geometry of a nondegenerate hermitian form on an $\Bbb R$- or $\Bbb C$-vector space $V$
from the perspective of the grassmannian $\Gr(k,V)$. The grassmannian consists of two complementary
parts: $\Gr^0(k,V)$, formed by nondegenerate subspaces, and the absolute.

\smallskip

The generic part $\Gr^0(k,V)$ of the grassmannian has several components, each corresponding to the
$k$-dimensional subspaces of $V$ of a certain nondegenerate signature. A tangent vector at a point
$p\in\Gr^0(k,V)$ is a linear map $t\in\Lin(p,p^\perp)$ and the trace $\tr(t_1^*t_2)$ provides an
hermitian metric on $\Gr^0(k,V)$. In particular, each component of $\Gr^0(k,V)$ is endowed with the
pseudo-riemannian metric $\Re\tr(t_1^*t_2)$. There exists, however, a more adequate way to
deal with the geometry of grassmannians. It involves {\it observations\/} and the {\it product.}

\smallskip

Extending by zero, the tangent vector $t\in\Lin(p,p^\perp)$ at $p\in\Gr^0(k,V)$ can be viewed as a
linear map $t\in\Lin(V,V)$. Conversely, a linear map $t\in\Lin(V,V)$ gives rise to the tangent vector
$$t_p:=\pi[p]t\pi'[p]$$
at $p\in\Gr^0(k,V)$, where $\pi'[p]:V\to p$ and $\pi[p]:V\to p^\perp$ stand for the orthogonal
projections onto $p$ and~$p^\perp$, respectively. The tangent vector $t_p$ is called the {\it
observation\/} of $t$ at $p$. Observations produce smooth vector fields on $\Gr^0(k,V)$ that play, in
the context of grassmannian geometries, a role similar to that played by left invariant vector fields
in the theory of Lie groups. Such vector fields are already important in the projective case $k=1$
[AGr1].

In the expression $\tr(t_1^*t_2)$ of the hermitian metric on $\Gr^0(k,V)$, one can take tangent vectors
$t_1,t_2$ observed at distinct points. Moreover, the characteristic polynomial $\ch(t_1^*t_2)$, or
simply the {\it product\/} $t_1^*t_2$ itself, can be taken instead of the hermitian metric. The
coefficients of $\ch(t_1^*t_2)$ provide not only the hermitian metric but also other geometric
characteristics such as the invariant $\det(t^*t)/\tr^k(t^*t)$ of the geodesic in $\Gr^0(k,V)$
determined by a tangent vector $t$ [AGr2]. Not infrequently, a fact involving the hermitian/pseudo-riemannian metric is just a
simple consequence of some property of the product. For example, the $m$-Pl\"ucker embedding
$$E^m:{\Gr}^0(k,V)\to{\Gr}^0\left({k\choose m},{\bigwedge}^mV\right),\quad p\mapsto{\bigwedge}^mp$$
is a (minimal) isometric embedding [AGr2]; this a trivial consequence of the functorial behaviour of
$E^m$ with respect to the product: $(E^mt_1)^*E^mt_2=E^m(t_1^*t_2)$. In other words, $E^m$ is an
`isometric embedding' in the sense of the product.

\smallskip

The variations involving observations and the product are endless. In order to measure distance in a
riemannian component $C$ of $\Gr^0(k,V)$, for example, we do not really need the metric, only
observations. Indeed, let $p_1,p_2\in C$ and let $t\in\Lin(p_1,p_1^\perp)\subset\Lin(V,V)$ be a nonnull
tangent vector at $p_1$. We observe $t$ at $p_2$ and then observe the result back at $p_1$ thus
obtaining a new tangent vector $t'$ at $p_1$. The change suffered by $t'$ reflects how $p_1$, $p_2$,
and $t$ are related and allows us to infer the distance $\dist(p_1,p_2)$. If, say, $k=1$ and $t$ is
tangent to the projective line joining $p_1$ and $p_2$, then $t'=\ta^2(p_1,p_2)t$, where
$$\ta(p_1,p_2):=\frac{\langle p_1,p_2\rangle\langle p_2,p_1\rangle}
{\langle p_1,p_1\rangle\langle p_2,p_2\rangle}$$
and $\dist(p_1,p_2)=\arccosh\sqrt{\ta(p_1,p_2}$ when the form on the subspace
$\Bbb Rp_1+\Bbb Rp_2\subset V$ is indefinite or $\dist(p_1,p_2)=\arccos\sqrt{\ta(p_1,p_2}$ when the
mentioned form is definite. Analogously, given points $p_1,p_2,p_3\in C$ and a tangent vector $t$ at
$p_1$, a similar procedure of successive observations leads to a new tangent vector $t'$ at $p_1$.
Again, $t'$ reflects the relation between $p_1,p_2,p_3$ and $t$. If $p_1,p_2,p_3\in C$ lie in a complex
projective line, i.e., the triangle $\Delta(p_1,p_2,p_3)$ is $\Bbb C$-plane, and $t$ is tangent to
this line, then $t'=t_{12}^2t_{23}^2t_{31}^2\exp\big(2i\Area\Delta(p_1,p_2,p_3)\big)t$, where
$t_{ij}:=\ta(p_i,p_j)$.

The above process of successive observations may be performed even when the metric is not riemannian
or when the points do not belong to a same component of $\Gr^0(k,V)$. In fact, in the global picture
composed of all the pseudo-riemannian pieces of $\Gr^0(k,V)$, observations and the product tie
everything together. Tangent vectors to a point in one piece are observable from the points in the
other pieces and we can take the product of tangent vectors at points in different pieces.

\smallskip

But what about the geometry on the absolute?

\smallskip

Our main objective in this paper is to show that observations and the product partially survive at a
degenerate point $p\in\Gr(k,V)$. Roughly speaking, let $q$ be the kernel of the hermitian form on $p$
and take the quotient $V_q:=q^\perp/q$. The points on the absolute of $\Gr(k,V)$ with the same $q$ form
a fibre of a certain bundle and observations and the product are defined for tangent vectors to such
fibre because $V_q$ is naturally equipped with a nondegenerate hermitian form. In particular, we obtain
a hermitian metric on the fibre.

When we consider the generic part of the absolute of $\Gr(k,V)$, that is, the points $p\in\Gr(k,V)$
such that the kernel of the hermitian form on $p$ is one-dimensional, the bundle in question fibres
over the absolute of a {\it projective\/} geometry $\Bbb PV$. Surprisingly, in this case, the hermitian
metric on the fibres provides the conformal (or conformal contact) structure on the absolute of
$\Bbb PV$. In other words, the conformal structure is exactly what remains from the metric when we
arrive at the absolute.

\smallskip

The basic geometrical objects in projective geometries, such as geodesics, totally geodesic subspaces,
equidistant loci, etc., have linear nature (see [AGr1] and references therein). In other words,
grassmannians constitute a natural habitat for important geometrical objects in classical spaces. Now,
we can see that the geometry of grassmannians also accounts for the relation between the
(pseudo-)riemannian metric of a projective geometry and the conformal structure on its absolute.
Particular cases of conformal structures that are amenable to this approach include the usual conformal
structure on the absolute of real hyperbolic space, the usual conformal structure on the absolute of de
Sitter space, the conformal contact structure on the absolute of complex hyperbolic space, and the
causal structure on the absolute of anti-de Sitter space (see Section 3).

The causal structure on the absolute of anti-de Sitter space plays an important role in general
relativity and the theory of black holes (see, for instance, [Car], [HoP], and references therein).
Another known application is possibly of a more unexpected nature: the absolute of anti-de Sitter
$4$-space can be identified with the lagrangian grassmannian $\Lambda(2)$, the space of lagrangian
planes in a $4$-dimensional symplectic $\Bbb R$-linear space. In other words, $\Lambda(2)$ is naturally
equipped with an extra structure arising from causal geometry (see [Cal] and Subsection 3.7). This
phenomenon may be seen as a geometrical manifestation of the exceptional isomorphism
$\Spin(2,3)=\Sp(4,\Bbb R)$.

\medskip

{\bf Acknowledgements.} It was with great sadness that we heard of the passing of Professor Waldyr
Rodrigues Jr. He was a dear friend and an exceptional scientist with whom we had the pleasure and the
privilege to discuss many subjects related to mathematics and physics. Waldyr was very fond of
natural and coordinate-free methods in geometry, a point of view that gave rise to this paper.

We are very grateful to Alexei L.~Gorodentsev and Nikolay A.~Tyurin for their stimulating interest in
this work and to the referees whose suggestions have greatly contributed to the exposition.

This paper was partially developed while the first author was enjoying the hospitality of the IHES
and, the third author, that of the MPIM.

\bigskip

\centerline{\bf2.~Stratification and fibre bundle}

\medskip

Let $V$ be an $n$-dimensional $\Bbb K$-vector space equipped with a
nondegenerate hermitian form $\langle-,-\rangle$ of arbitrary signature, where $\Bbb K=\Bbb R$
or $\Bbb K=\Bbb C$. The grassmannian $\Gr_\Bbb K(k,V)$ of
$k$-dimensional $\Bbb K$-vector subspaces in $V$ can be described as
follows. Fix a $\Bbb K$-vector space $P$ such that
$\dim_\Bbb KP=k$. Denote by
$$M:=\big\{p\in{\Lin}_\Bbb K(P,V)\mid\ker p=0\big\}$$
the open subset of all monomorphisms in the $\Bbb K$-vector space
$\Lin_\Bbb K(P,V)$. The group $\GL_\Bbb KP$ acts from the right on
$\Lin_\Bbb K(P,V)$ and on $M$. The grassmannian $\Gr_\Bbb K(k,V)$ is
simply the quotient space
$${\Gr}_\Bbb K(k,V):=M/{\GL}_\Bbb KP,\qquad\pi:M\to M/{\GL}_\Bbb KP.$$

We will not distinguish between the notation of points in $\Gr_\Bbb K(k,V)$ and of their
representatives in~$M$. Moreover, we will frequently write $p$ in place of the image $pP$ and
$p^\perp$, in place of the orthogonal $(pP)^\perp$ to that image. For example, $V/p$ will denote
$V/pP$.

\smallskip

The tangent space $\T_pM$ is usually identified with
$\Lin_\Bbb K(P,V)$ as follows. Given $\varphi:P\to V$, the curve
$c(\varepsilon):=p+\varepsilon\varphi$ lives in $M$ for small
$\varepsilon$ and $\varphi$ is identified with $\dot c(0)$. We use a
slightly different identification $\T_pM=\Lin_\Bbb K(p,V)$ where $p$
takes the place of $P$. In this way, $\overline t\in\Lin_\Bbb K(p,V)$
is interpreted as the tangent vector $\dot c(0)$, where
$c(\varepsilon):=(1+\varepsilon\widetilde t)p$ and
$\widetilde t\in\Lin_\Bbb K(V,V)$ extends $\overline t$. Note that
$\overline t\in\Lin_\Bbb K(p,V)$ is tangent at $p$ to the orbit
$p\GL_\Bbb KP$ if and only if $\overline tp\subset p$ since in this
case $(1+\varepsilon\widetilde t)p=pg$ for small $\varepsilon$ and
suitable $g\in\GL_\Bbb KP$. Therefore,
$\T_p\Gr_\Bbb K(k,V)=\Lin_\Bbb K(p,V/p)$.

\medskip

{\bf2.1.~Remark.} Let $w\subset V$ be a linear subspace and let
$\Gr_\Bbb K(k,w,V)\subset\Gr_\Bbb K(k,V)$ stand for the space of all $k$-dimensional linear subspaces
in $V$ that are included in $w$. Then ${\T}_p{\Gr}_\Bbb K(k,w,V)={\Lin}_\Bbb K(p,w/p)$ for all
$p\in\Gr_\Bbb K(k,w,V)$.

Dually, let $\Gr_\Bbb K(k,V,q)\subset\Gr_\Bbb K(k,V)$ denote the space of all $k$-dimensional linear
subspaces in $V$ containing $q$, where $q\subset V$ is a given $d$-dimensional subspace. Then
${\T}_p{\Gr}_\Bbb K(k,V,q)={\Lin}_\Bbb K(p/q,V/p)$ for all $p\in\Gr_\Bbb K(k,V,q)$.
\hfill$_\blacksquare$

\medskip

The degree of degeneracy of a point $p\in\Gr_\Bbb K(k,V)$ is the dimension of the kernel of the
hermitian form restricted to $p$, that is, the dimension of the subspace $p\cap p^\perp\subset V$. The
grassmannian $\Gr_\Bbb K(k,V)$ is stratified according to such degree:
$${\Gr}_\Bbb K(k,V)=\bigsqcup\limits_d{\Gr}_\Bbb K^d(k,V),\qquad
{\Gr}_\Bbb K^d(k,V):=\big\{p\in{\Gr}_\Bbb K(k,V)\mid\dim_\Bbb K(p\cap
p^\perp)=d\big\}.$$
Subspaces of $V$ of a given signature form an $\U V$-orbit. Therefore, each stratum is the disjoint
union of a finite number of such orbits, hence, a manifold (note that a stratum is not necessarily
connected, so we use the word `stratification' with a meaning slightly different from the usual one).
There are finitely many strata and the closure of the stratum $\Gr_\Bbb K^{d_0}(k,V)$ is the disjoint
union of strata $\bigsqcup\limits_{d\geqslant d_0}\Gr_\Bbb K^d(k,V)$. The generic part
$\Gr_\Bbb K^0(k,V)$ corresponding to the nondegenerate subspaces of $V$ is open in $\Gr_\Bbb K(k,V)$.

\smallskip

Associating to each $p\in\Gr_\Bbb K(k,V)$ the kernel of the hermitian form on $p$, we get the
$\U V$-equivariant fibre bundle
$$\pi_d:{\Gr}_\Bbb K^d(k,V)\to{\Gr}_\Bbb K^d(d,V),\qquad\pi_d:p\mapsto p\cap p^\perp.\leqno{(2.2)}$$
The fibre $\pi_d^{-1}(q)$ can be naturally identified with $\Gr_\Bbb K^0(k-d,V_q)$, where
$V_q:=q^\perp/q$ is equipped with a natural nondegenerate hermitian form and $\dim_\Bbb KV_q=n-2d$.

\smallskip

We are particularly interested in the case $d=1$ corresponding to the generic part of the absolute of
$\Gr_\Bbb K(k,V)$. This is the case that will be applied, in the next section, to the study of
conformal structures on the absolute of {\it projective\/} geometries: when $d=1$, the
bundle (2.2) fibres over $\Gr_\Bbb K^1(1,V)$ which is the absolute, denoted by $\S V$, of
$\Bbb P_\Bbb KV=\Gr_\Bbb K(1,V)$. Fibres are easy to visualize, as they correspond to subspaces
`rotating' about their common one-dimensional kernel, hence, forming a dense open part of the
grassmannian $\Gr_\Bbb K(k-1,n-2)$. The bundle itself is also simple to describe, as the next couple of
propositions show.

\medskip

{\bf2.3.~Proposition.} {\it Let\/ $\Bbb K=\Bbb R$. The bundle\/
$\pi_1:\Gr_\Bbb R^1(k,V)\to\Gr_\Bbb R^1(1,V)$ is the nondegenerate part of a grassmannization of the
tangent bundle of the absolute\/ $\S V$ of the projective geometry\/ $\Bbb P_\Bbb RV$.}

\medskip

{\bf Proof.} Take $q\in\Gr_\Bbb R^1(1,V)=\S V$ and define $V_q:=q^\perp/q$. Since the fibre
$\pi_1^{-1}(q)$ is naturally identified with the nondegenerate part $\Gr_\Bbb R^0(k-1,V_q)$ of the
grassmannian $\Gr_\Bbb R(k-1,V_q)$ and the subspace $q\subset V$ is one-dimensional, it suffices to
prove that the tangent space to $\S V$ at $q$ has the form $\T_q\S V=\Lin_\Bbb R(q,V_q)$.

Take $u\in V$ such that $\S V\subset\Bbb P_\Bbb RV$ is locally given by the equation $f(x)=0$ in a
neighbourhood of $q$, where
$f(x):=\displaystyle\frac{\langle x,x\rangle}{\langle x,u\rangle\langle
u,x\rangle}$.
Let $t\in\Lin(q,V/q)$ be a tangent vector to $\Bbb P_\Bbb RV$ at $q$ and let $\widetilde t\in\Lin(V,V)$
be a lift of $t$. Then $t$ is tangent to $\S V$ if and only if $tf=0$, i.e.,
$$\frac d{d\varepsilon}\Big|_{\varepsilon=0}\frac{\langle
q+\varepsilon\widetilde tq,q+\varepsilon\widetilde tq\rangle}{\langle
q+\varepsilon\widetilde tq,u\rangle\langle u,q+\varepsilon\widetilde
tq\rangle}=\frac{2\Re\langle\widetilde tq,q\rangle}{\langle
q,u\rangle\langle u,q\rangle}=\frac{2\Re\langle tq,q\rangle}{\langle
q,u\rangle\langle u,q\rangle}=0.$$
(In the formula, $q\in V$ denotes an element representing the subspace
$q$.)
\hfill$_\blacksquare$

\medskip

When $\Bbb K=\Bbb C$, a calculus analogous to the one presented in the proof of Proposition 2.3 shows
that the tangent space to the absolute $\S V$ of the projective geometry
$\Bbb P_\Bbb CV=\Gr_\Bbb C(1,V)$ at a point $q\in\S V$ is given by
$${\T}_q{\S}V=\big\{t\in{\Lin}_\Bbb C(q,V/q)\mid\Re\langle tq,q\rangle=0\big\}.$$
Hence, $\Lin_\Bbb C(q,V_q)\subset\T_q\S V$ is not the entire tangent space at $q$ as in the real case,
but a maximal complex subspace of $\T_q\S V$. In other words, the subspaces of the form
$\Lin_\Bbb C(q,V_q)$ give rise to a CR-distribution on $\S V$. Summarizing, we have the following
proposition.

\medskip

{\bf2.4.~Proposition.} {\it Let\/ $\Bbb K=\Bbb C$. The bundle\/
$\pi_1:\Gr_\Bbb C^1(k,V)\to\Gr_\Bbb C^1(1,V)$ is the nondegenerate part of a grassmannization of the\/
{\rm CR-}distribution in\/ $\T\S V$ given by\/ $\Lin_\Bbb C(q,V_q)\subset\T_q\S V$, where\/ $q\in\S V$
and\/ $V_q:=q^\perp/q$.}

\medskip

The case $d>1$ deals with a `degenerate' part of the absolute of the
grassmannian. The bundle $\pi_d:\Gr_\Bbb K^d(k,V)\to\Gr_\Bbb K^d(d,V)$
is no longer the grassmannization of the tangent bundle of
$\Gr_\Bbb K^d(d,V)$ (considering, say, the real case). Nevertheless,
the bundle indicates a distinguished distribution whose geometric
nature would be interesting to dwell on.

\medskip

The fact that $\pi_d^{-1}(q)=\Gr_\Bbb K^0(k-d,V_q)$, where $V_q:=q^\perp/q$ is taken with its natural
induced nondegenerate hermitian form, means that the fibres of the bundle (2.2) are naturally
equipped with the geometric structure consisting of observations and of the product (this structure on
the nondegenerate part of a grassmannian was discussed in the introduction).

Observations and the product on the fibre $\pi_d^{-1}(q)$ can be explicitly described as follows. By
Remark~2.1,
$${\T}_p\pi_d^{-1}(q)={\Lin}_{\Bbb K}(p/q,q^\perp/p)$$
for all $p\in\pi_d^{-1}(q)$ because $\Gr_\Bbb K^0(k-d,V_q)$ is open in $\Gr_\Bbb K(k-d,V_q)$. Let
$p\in\pi_d^{-1}(q)$. Clearly, $q=p\cap p^\perp$, $q^\perp=p+p^\perp$, and $V_q=p_0\oplus p_0^\perp$,
where $p_0:=p/q$ and $p_0^\perp=p^\perp/q$. We have $\T_p\pi_d^{-1}(q)=\Lin_\Bbb K(p_0,p_0^\perp)$.
Denote by $\pi'[p_0]$ and $\pi[p_0]$ the orthogonal projections onto $p_0$ and $p_0^\perp$,
respectively. Then the linear map $t\in\Lin_\Bbb K(V_q,V_q)$ gives rise to the tangent vector
$$t_{p_0}:=\pi[p_0]t\pi'[p_0]$$
to $\pi_d^{-1}(q)$ at $p_0$ (the observation of $t$ at $p_0$). The product $t_1^*t_2$, where
$t_1,t_2\in\Lin_\Bbb K(V_q,V_q)$, gives rise to the hermitian metric $\tr(t_1^*t_2)$ on the fibre when
$t_1$ and $t_2$ are observed from a same point. In particular, the fibre $\pi_d^{-1}(q)$ consists of a
collection of pseudo-Riemannian manifolds, each corresponding to the $(k-d)$-dimensional subspaces of
$V_q$ of a given signature. Some particular examples will be considered in the next section.

\bigskip

\centerline{\bf3.~Conformal structures}

\medskip

The (pseudo-)riemannian metric on a space typically corresponds to some geometric structure on its
ideal boundary (absolute). The correspondence is such that both the metric and the boundary structure
have the same automorphisms. For example, the absolute of real hyperbolic space is a conformal sphere,
the absolute of de Sitter space is a conformal sphere, the absolute of complex hyperbolic space is a
CR-sphere, and the absolute of anti-de Sitter space has a causal structure.

We will see in this section that the bundle (2.2) allows us to explain these `gravity/conformal structure'
correspondences. At the end of the day, the hermitian form on $V$ is responsible not only for the
metric, but also for the geometric structure on the absolute. The role of the bundle (2.2) is to
provide a way for this hermitian form to `reach' the tangent space of the absolute.

\smallskip

We call the bundle $\pi_1:\Gr_\Bbb K^1(k,V)\to\Gr_\Bbb K^1(1,V)$ the {\it conformal\/} ({\it conformal
contact\/} when $\Bbb K=\Bbb C$) structure on the absolute $\S V=\Gr_\Bbb K^1(1,V)$ of the projective
geometry $\Bbb P_\Bbb KV=\Gr_\Bbb K(1,V)$. Let us analyze the conformal (conformal contact) structure
for $k=2$ in a few cases.

\smallskip

We will use the following observation in all the examples below. Let
$p\in\Gr_\Bbb K^0(1,V)$ and let $t_1,t_2\in\Lin_\Bbb K(p,p^\perp)$ be tangent vectors at
$p$. It is easy to see that, in these projective settings, the expression
$$\pm\tr(t_1^*t_2)=\pm\frac{\langle t_1p,t_2p\rangle}{\langle p,p\rangle}\leqno(3.1)$$
gives the (pseudo-)riemannian metrics in the components of $\Gr_\Bbb K^0(1,V)$.

\smallskip

{\bf3.2.~Real hyperbolic and de Sitter spaces.}~Let $\Bbb K=\Bbb R$ and assume that $V$ is equipped
with an hermitian (i.e., bilinear symmetric) form of signature $+\dots+-$. Let $n:=\dim_\Bbb RV$. Then $\Gr_\Bbb R^0(1,n)$
has two components: one, called the {\it real hyperbolic space\/} $\Bbb H_\Bbb R^{n-1}$, consists of
negative points; the other, called the {\it de Sitter space\/} $\d\Bbb S^{n-1}$, consists of positive
points. The absolute equals $\Gr_\Bbb R^1(1,n)$ and is the sphere $\Bbb S^{n-2}$ of null points. In
other words, the real hyperbolic space and the de Sitter space are glued along their absolutes; we call
the projective space $\Gr_\Bbb R(1,V)$ the {\it extended\/} real hyperbolic space.

Let $p\in\Gr_\Bbb R^0(1,V)$ and let $t_1,t_2\in\Lin_\Bbb R(p,p^\perp)$ be tangent vectors at $p$. We
introduce the metric $\langle t_1,t_2\rangle:=-\tr(t_1^*t_2)$ on $\Gr_\Bbb R^0(1,V)$. Due to (3.1) and to the
signature of the form, this is a riemannian metric on $\Bbb H_\Bbb R^{n-1}$ and a lorentzian metric on
$\d\Bbb S^{n-1}$. In this way, we obtain the usual constant curvature metrics on the hyperbolic and de
Sitter spaces [AGr1].

Let $q\in\S V$. The fibre $\pi_1^{-1}(q)$ of the conformal structure
$\pi_1:\Gr_\Bbb R^1(2,V)\to\Gr_\Bbb R^1(1,V)$ is naturally identified with
$\Gr^0_\Bbb R(1,V_q)=\Bbb P_\Bbb R^{n-3}$ because $k=2$, $d=1$, and $\dim_\Bbb RV_q=n-2$,
where $V_q:=q^\perp/q$. The
induced form on $V_q$ is positive-definite since the signature of the form on $q^\perp$ is $0+\dots+$.
The form on $V_q$ therefore provides a riemannian metric on each fibre. The latter is a metric of
constant curvature on~$\Bbb P_\Bbb R^{n-3}$~[AGr1].

By Proposition 2.3, the bundle $\pi_1:\Gr_\Bbb R^1(2,V)\to\Gr_\Bbb R^1(1,V)$ is the projectivization
of the tangent bundle of the absolute $\Bbb S^{n-2}$. Distances in the fibre are nothing but the angles
of the standard\footnote{If one wishes to deal with angles varying in $[0,2\pi]$, then the
projectivization should be taken with respect to $\Bbb R^+$, from the very beginning.}
conformal structure on the sphere $\Bbb S^{n-2}$.
\hfill$_\blacksquare$

\medskip

{\bf3.3.~Anti-de Sitter space.}~Take $\Bbb K=\Bbb R$ and $V$ of signature $+\dots+--$. The {\it
anti-de Sitter space\/} $\ad\Bbb S^{n-1}$ is the negative part of $\Bbb P_\Bbb RV$. It is a lorentzian
manifold. The fibre $\pi_1^{-1}(q)$ of the conformal structure
$\pi_1:\Gr_\Bbb R^1(2,V)\to\Gr_\Bbb R^1(1,V)$ is naturally identified with $\Gr^0_\Bbb R(1,V_q)$ and
is an open dense set in $\Bbb P_\Bbb R^{n-3}$. Since the signature of the form on $q^\perp$ is
$0-+\dots+$, the signature of the form on $V_q=q^\perp/q$ is
$-+\dots+$ and the fibres carry the metric of the extended hyperbolic space of the previous example
(each fibre is an extended hyperbolic space minus its absolute).

Such a conformal structure is known as a {\it causal\/} structure. It plays in anti-de Sitter geometry
the same role played by the standard conformal strucure in the absolute of the real hyperbolic space.

Curiously, the absolute of anti-de Sitter $4$-space $\ad\Bbb S^4$ can be
naturally identified with the space $\Lambda(2)$ of lagrangian planes in a
$4$-dimensional symplectic $\Bbb R$-linear space (see Subsection 3.6). Hence,
the {\it lagrangian grassmannian\/} $\Lambda(2)$ possesses a natural causal
structure.
\hfill$_\blacksquare$

\medskip

For $k>2$, the geometry on the grassmannization
$\Gr_\Bbb R^1(k,V)\to\Gr_\Bbb R^1(1,V)$ of the tangent bundle of the
absolute (see Proposition 2.3) is related to the case $k=2$ in the same
way as are related the grassmannian and the projective
geometries.

\medskip

{\bf3.4.~Complex hyperbolic space.} Take $\Bbb K=\Bbb C$ and $V$ of signature $+\dots+-$. The {\it
complex hyperbolic space\/} $\Bbb H_\Bbb C^{n-1}$ is the negative part of $\Bbb P_\Bbb CV$ and its ideal boundary,
the absolute, is the sphere $\S V\simeq\Bbb S^{2n-3}$. Every fibre
$\Bbb P_\Bbb RV_q\simeq\Bbb P_\Bbb C^{n-3}$ of the conformal contact structure carries the Fubini-Study
metric. By Proposition 2.4, the conformal contact structure $\pi_1:\Gr_\Bbb C^1(2,V)\to\Gr_\Bbb C^1(1,V)$ is the
projectivization of a CR-distribution in the tangent bundle of the absolute $\Bbb P_\Bbb C^{n-3}$ and
the distances in a fibre are the angles between complex directions.
\hfill$_\blacksquare$

\medskip

Note that any projective geometry can play the role of a
(contact) conformal structure. In particular, the (contact) conformal structure can possess
its own absolute and so on $\dots$

\medskip

{\bf3.5.~Comments and questions.} Algebraic formulae that deal with
geometrical quantities  (the {\it tance\/} and the $\eta${\it-invariant} are some simple
examples, see [AGr1, Example 5.12]) work well for points in distinct
pieces of $\Gr_\Bbb K^0(k,V)$. Such formulae use
to alter its geometrical sense when the points involved are taken in different components of
$\Gr_\Bbb K^0(k,V)$. In this respect, it is interesting to understand if
there is an explicit geometrical interpretation of observations and the product
in the terms of the usual (pseudo-)riemannian concepts for $k>2$.

The bundle $\pi_d:\Gr_\Bbb K^d(k,V)\to\Gr_\Bbb K^d(d,V)$ might admit a
canonical connection. If so, what is its explicit description?

\medskip

{\bf3.6.~The lagrangian grassmannian.} The material in this subsection
is well-known and closely related to that in [Cal]. Our intention is
simply to give it a complete description in the spirit of [AGr1] and of
the present paper.

Let $V$ be a $4$-dimensional $\Bbb R$-linear space. We fix
$0\ne\omega\in\bigwedge^4V$. The wedge product endows $\bigwedge^2V$
with the symmetric bilinear form
$$\langle-,-\rangle:{\bigwedge}^2V\times{\bigwedge}^2V\to\Bbb R,\qquad
a\wedge b=\langle a,b\rangle\omega.$$
The form $\langle-,-\rangle$ is nondegenerate of signature $---+++$.
Indeed, if $b_1,b_2,b_3,b_4$ is a basis in $V$ such that
$b_1\wedge b_2\wedge b_3\wedge b_4=\omega$, then the elements
$b_1\wedge b_2\pm b_3\wedge b_4$, $b_1\wedge b_3\pm b_2\wedge b_4$, and
$b_1\wedge b_4\pm b_2\wedge b_3$ form an orthogonal basis in
$\bigwedge^2V$ from which the signature of the form can be easily
inferred.

\smallskip

Given a $2$-dimensional subspace $W\leqslant V$, the restriction of
$\langle-,-\rangle$ to $\bigwedge^2W\leqslant\bigwedge^2V$ is null. In
other words, the Pl\"ucker map
$${\Gr}_\Bbb R(2,V)\to{\Gr}_\Bbb R\Big(1,{\bigwedge}^2V\Big)=\Bbb
P{\bigwedge}^2V$$
embeds $\Gr_\Bbb R(2,V)$ into the absolute $\S V$ of the projective
geometry $\Bbb P{\bigwedge}^2V$. Hence,
$\Gr_\Bbb R(2,V)\to\S V$ is a diffeomorphism, because
$\Gr_\Bbb R(2,V)$ and $\S V$ are manifolds of the same dimension, $\Gr_\Bbb R(2,V)$ is
compact without boundary, and $\S V$ is connected (alternatively, one
may simply note that the conditions $0\ne a\in{\bigwedge}^2V$ and
$a\wedge a=0$ imply  $a=v\wedge w$ for some linearly independent
$v,w\in V$).

\smallskip

Let $0\ne\eta\in{\bigwedge}^2V$ be negative, that is,
$\langle\eta,\eta\rangle<0$. Then $V$ is equipped with the symplectic form
$$(-,-):V\times V\to\Bbb R,\qquad (v,w):=\langle\eta,v\wedge w\rangle.$$
A $2$-dimensional subspace $W\leqslant V$ is called {\it lagrangian\/} if
the restriction of $(-,-)$ to $W\times W$ is null, that is, if
$(-,-)\big\vert_{W\times W}=0$. The subset $\Lambda(2)\subset\Gr_\Bbb R(2,V)=\S V$ corresponding
to lagrangian planes is known as the {\it lagrangian grassmannian.} Clearly,
$\Lambda(2)=\Bbb P\eta^\perp\cap\S V$, where $\eta^\perp$ stands for
the orthogonal complement to $\eta$ with respect to the symmetric bilinear
form $\langle-,-\rangle$. It follows from the description of the tangent space to a point in $\S V$
in the proof of Proposition 2.3 and [AGG, Lemma~4.2.2] that
the intersection of the hypersurfaces $\Bbb P\eta^\perp$ and $\S V$ in
$\Bbb P\bigwedge^2V$ is transversal. So, the lagrangian grassmannian is a
$3$-dimensional submanifold of $\Gr_\Bbb R(2,V)$.

The symmetric bilinear form restricted to $\eta^\perp$ has signature $--+++$ and
$\Bbb P\eta^\perp\cap\S V$ is nothing but the absolute of the projective
geometry $\Bbb P\eta^\perp$. In other words, the lagrangian grassmannian can be
naturally identified with the absolute of anti-de Sitter $4$-space. It is therefore
endowed with a causal structure.

Causal structures in groups of symplectomorphisms/contactomorphisms constitute an important
topic in the geometry of contact manifolds [ElP].

\bigskip

\centerline{\bf4.~References}

\medskip

[AGG] S.~Anan$'$in, C.~H.~Grossi, N.~Gusevskii, {\it Complex hyperbolic structures
on disc bundles over surfaces,} International Mathematics Research Notices, 2011
{\bf19}, 4295--4375

\smallskip

[AGr1] S.~Anan$'$in, C.~H.~Grossi, {\it Coordinate-free classic geometries,}
Moscow Mathematical Journal, 2011 {\bf11}, 633--655

\smallskip

[AGr2] S.~Anan$'$in, C.~H.~Grossi, {\it Differential geometry of grassmannians
and Pl\"ucker map,} Central European Journal of Mathematics, 2012 {\bf3}, 873--884

\smallskip

[Cal] D.~Calegari, {\it Causal geometry,} Geometry and the imagination (blog post
at lamington.\break wordpress.com/2009/12/10/causal-geometry), 2009

\smallskip

[Car] S.~Carlip, {\it The\/ $(2+1)$-dimensional black hole,} Classical and
quantum gravity, 1995 {\bf12}, 2853--2880

\smallskip

[ElP] Y.~Eliashberg, L.~Polterovich, {\it Partially ordered groups and
geometry of contact transformations,} Geometric and Functional Analysis
(GAFA), 2000 {\bf10}, 1448--1476

\smallskip

[HoP] S.~Holst, P.~Peld\'an, {\it Black holes and causal structure in anti-de Sitter
isometric spacetimes,} Classical and quantum gravity, 1997 {\bf14}, 3433--3452

\enddocument